\documentclass[twoside]{article}
\usepackage{amsmath}
\usepackage{url}
\usepackage{indentfirst}
\usepackage{fancyhdr}
\usepackage[margin=4.76cm]{geometry}
\usepackage{graphicx}
\usepackage{subfigure}
\begin{document}
\pagestyle{fancy}
\renewcommand{\headrulewidth}{0.4pt}
\fancyhead{}
\fancyhead[CE]{\textsc{\footnotesize Meng Zhang}}
\fancyhead[CO]{\textsc{\footnotesize The Fractal and The Recurrence Equations Concerning The Integer Partitions}}
\title{The Fractal and The Recurrence Equations Concerning The Integer Partitions}
\author{Meng Zhang}
\date{}
\maketitle

\begin{abstract}
This paper introduced a way of fractal to solve the problem of taking count of the integer partitions, furthermore, using the method in this paper some recurrence equations concerning the integer partitions can be deduced, including the pentagonal number theorem.
\end{abstract}
\setlength{\parskip}{0.5\baselineskip}
\begin{small}
\end{small}
\vspace{1.5ex}
\section{Introduction}

\noindent A partition of a positive integer $\displaystyle{n}$ is a way of writing $\displaystyle{n}$ as a sum of positive integers. The number of partitions of $\displaystyle{n}$ is given by the partition function $\displaystyle{p(n)}$. In the history of researching the expression of $\displaystyle{p(n)}$, the methods such as generating functions and complex analysis are widely used. The fractal can not only describe lots of the natural phenomena, but also work in the recurrence equations. In this paper a new method of fractal will be used to deduce some recurrence equations concerning the integer partitions, including the pentagonal number theorem.

\section{Preliminaries}

\noindent It is well known that the Mandelbrot set \textsuperscript{\cite{1}} is generated by the recurrence equation $Z_{n+1}$ = $Z_n^2$ + C, where $Z_n$ is a complex number and C is a complex constant \textsuperscript{\cite{1}\cite{2}}.In this paper let us call the recurrence equation $Z_{n+1}$ = $Z_n^2$ + C \textsuperscript{\cite{2}} the generator of the Mandelbrot set. Next we should find the generator of the integer partition functions. It is defined that $\displaystyle{p(m)=0}$ if $\displaystyle{m}$ $<$ 0 and $\displaystyle{p(0)=1}$. Let us call the value $\displaystyle{m}$ in $\displaystyle{p(m)}$ the tab of $\displaystyle{p(m)}$. Here let us write an equation below,\vspace{1.5ex}
\\ $\displaystyle{p(n)=\{p(0)\}+\{p(1)\}+\{p(2)\}}$ + $\cdots$ + $\displaystyle{\{p(i)\}}$ + $\cdots$ + $\displaystyle{\{p(n-1)\}}$,\vspace{1.5ex}\hfill(1)
\\where $\displaystyle{\{p(i)\}}$ is called a parcel here.

To accomplish the proof of the main theorem in the next section, some conceptions are required. Let us call the elements in the symbol \{ \} and the symbol itself a parcel. Every parcel may generate other parcel(s). Let us call a parcel which has generated its child-parcel(s) a father-parcel. A parcel has a father-parcel if it is generated by other parcel.

Let us define that a parcel is surrounded with the symbol \{ \} if it has not generated its child-parcel(s), otherwise the parcel and its child-parcel(s) are surrounded with the symbol [ ], let us call the elements in the symbol [ ] and the symbol itself a cell. The symbol [ ] has the same meaning with the brackets in arithmetic. In this paper the symbol [ ] is different from the symbol $\lfloor $ $\rfloor$, $\lfloor x \rfloor$ denotes a rounding down function which will be used later. The child-parcel(s) are the repeated count in their father-parcel, they should be subtracted from their father-parcel. Here let us define that a parcel will lost the symbol \{ \} out of itself if it has generated its child-parcel(s). Next we should find the expression of every parcel.

\section{Deducing the Generator}

\noindent \emph{Definition 1}

\emph{$\displaystyle{\Omega(n)}$ is a function that can let n be written as an arbitrary partition.}

For example, $\displaystyle{\Omega(15)=7+8}$ is an expression of the partitions of integer 15. The number of representations of $\displaystyle{\Omega(n)}$ is $\displaystyle{p(n)}$.

\noindent \emph{Definition 2}

\emph{$\displaystyle{m}$ is called a fixed number if $\displaystyle{n = m+\Omega(n-m)}$, where $\displaystyle{m}$ can not be divided.}

Let us assume the sum of the child-parcel(s) that a parcel will generate is $\displaystyle{\sum_{i=0}^{\lambda} \{p(i)\}}$, we should find the $\displaystyle{\lambda}$ in it.\vspace{2ex}

\noindent \emph{The Main Theorem}

\emph{In a cell, let $\displaystyle{\lambda}$ denote the quantity of the child-parcel(s) that a
parcel $\displaystyle{\{p(\tau)\}}$ will generate, let $\displaystyle{n}$ denote the tab in its father-parcel, then $\displaystyle{\lambda=2\tau-n}$. $\displaystyle{\lambda=0}$ if $\displaystyle{2\tau-n\leq0}$.} \vspace{2ex}

\noindent Proof

When $\displaystyle{n}$ is an odd number, let $\displaystyle{t}$ and $\displaystyle{\{p(n-t)\}}$ be bijection. When 1 $\leq$ $\displaystyle{t}$ $<$ $\displaystyle{n/2}$,
let $\displaystyle{t}$ be a fixed number, the count of partition mapped $\displaystyle{t}$ is $\displaystyle{p(n-t)}$, here let us sign it \{$\displaystyle{p(n-t)}$\} because $\displaystyle{p(n-t)}$ has repeat count.

Let us divide $\displaystyle{n}$ into $\displaystyle{t}$ and $\displaystyle{n-t}$, let $\displaystyle{\xi(i)}$ denote the number in the section [$\displaystyle{t}$, $\displaystyle{n-t}$], where $\displaystyle{t}$ $\leq$ $\displaystyle{i}$ $\leq$ $\displaystyle{n-t}$. Because $\displaystyle{t+\xi(i)}$ $\leq$ $\displaystyle{n}$, let $\varphi$ denote the largest $\displaystyle{\xi(i)}$, we have
$\displaystyle{\varphi=n-t}$. 	                                                           \hfill(2)

$\varphi$ = $\displaystyle{n-t}$ $>$ $\displaystyle{n/2}$ $>$ $\displaystyle{t}$ because $\displaystyle{t}$ $<$ $\displaystyle{n/2}$. Now let us regard $\varphi$ as the fixed number, the partition count of $\displaystyle{n=\varphi+t+\Omega(n-t-\varphi)}$ can be mapped as
$\displaystyle{p(n-t-\varphi)}$ because it has been counted when $\varphi$ is a fixed number. Using formula (2) we have $\displaystyle{p(n-t-\varphi)=\{p(0)\}}$.

Also, we can regard every number in the section ($\displaystyle{t}$, $\varphi$) as a fixed number.

Therefore, in the section [$\displaystyle{t}$, $\varphi$], every number can be mapped by order as\vspace{2ex}
$\displaystyle{\{p(n-t-\varphi)\}=\{p(0)\}}$ \quad\quad\;\;\;\,$\leftrightarrow$\quad\quad\;\;\;\, $\displaystyle{n=\varphi+t+\Omega(n-t-\varphi)}$,
\\$\displaystyle{\{p(n-t-(\varphi-1))\}=\{p(1)\}}$ \quad\;\, $\leftrightarrow$ \quad\;\,  $\displaystyle{n=(\varphi-1)+t+\Omega(n-t-(\varphi-1))}$,
\\$\displaystyle{\{p(n-t-(\varphi-2))\}=\{p(2)\}}$ \quad\;\, $\leftrightarrow$ \quad\;\,   $\displaystyle{n=(\varphi-2)+t+\Omega(n-t-(\varphi-2))}$,
\\$\cdots$
\\$\displaystyle{\{p(n-t-(t+1))\}=\{p(n-2t-1)\}}$\;$\leftrightarrow$\;$\displaystyle{n=(t+1)+t+\Omega(n-t-(t+1))}$,\vspace{1ex}
\\where ``$\leftrightarrow$'' is a symbol of bijection.

Then let us calculate the count of the equations above, we have

$\displaystyle{\lambda=n-2t-1-0+1=n-2t}$.

Thus,

$\displaystyle{\lambda=n-2t=n-2(n-\tau)=2\tau-n}$ because $\displaystyle{\tau=n-t}$.

When $\displaystyle{n/2}$ $\leq$ $\displaystyle{t}$ $\leq$ $\displaystyle{n}$, there is no number in the section $\displaystyle{[t, n-t]}$, thus $\displaystyle{\lambda=0}$ because $\displaystyle{n-2t\leq0}$ if $\displaystyle{(2\tau-n)\leq0}$.

In the situation of even, the process of analysis is as same as odd. This concludes the proof of the main theorem.

For example, let us use the generator to calculate $\displaystyle{p(10)}$. Let us regard $\displaystyle{p(10)}$ as a father-parcel, using the main theorem (or it can be called the generator), at first, every parcel has not generated its child-parcel(s), so they are surrounded by the symbol \{ \}, next, every parcel generates its child-parcel(s) and has become a cell, the child-parcel(s) are arranged by order according to their tabs, their first child-parcel should be $\displaystyle{p(0)}$. The quantity of the child-parcel(s) that a parcel will generate is given by the main theorem. Therefore the tab of the last child-parcel equals the quantity of the child-parcel(s) minus 1. The process of generating will be continued until there is no parcel can generate its child-parcel(s). We can ``zoom in'' $\displaystyle{p(10)}$ entirely with 5 steps, which are shown below.\vspace{2ex}
\\$\displaystyle{p}$(10)=$\displaystyle{\sum_{i=0}^{9}}$$\displaystyle{\{p(i)\}}$
\\\hspace*{8mm}=$\displaystyle{p}$(0)+$\displaystyle{p}$(1)+$\displaystyle{p}$(2)+$\displaystyle{p}$(3)+$\displaystyle{p}$(4)+$\displaystyle{p}$(5)+[$\displaystyle{p}$(6)$\displaystyle{-}$\{$\displaystyle{p}$(0)\}$\displaystyle{-}$\{$\displaystyle{p}$(1)\}]+[$\displaystyle{p}$(7)$\displaystyle{-}$\{$\displaystyle{p}$(0)\}
\\\hspace*{10.5mm}$\displaystyle{-}$\{$\displaystyle{p}$(1)\}$\displaystyle{-}$\{$\displaystyle{p}$(2)\}$\displaystyle{-}$\{$\displaystyle{p}$(3)\}]+[$\displaystyle{p}$(8)$\displaystyle{-}$\{$\displaystyle{p}$(0)\}$\displaystyle{-}$\{$\displaystyle{p}$(1)\}$\displaystyle{-}$\{$\displaystyle{p}$(2)\}$\displaystyle{-}$\{$\displaystyle{p}$(3)\}
\\\hspace*{10.5mm}$\displaystyle{-}$\{$\displaystyle{p}$(4)\}$\displaystyle{-}$\{$\displaystyle{p}$(5)\}]+[$\displaystyle{p}$(9)$\displaystyle{-}$\{$\displaystyle{p}$(0)\}$\displaystyle{-}$\{$\displaystyle{p}$(1)\}$\displaystyle{-}$\{$\displaystyle{p}$(2)\}$\displaystyle{-}$\{$\displaystyle{p}$(3)\}$\displaystyle{-}$\{$\displaystyle{p}$(4)\}
\\\hspace*{10.5mm}$\displaystyle{-}$\{$\displaystyle{p}$(5)\}$\displaystyle{-}$\{$\displaystyle{p}$(6)\}$\displaystyle{-}$\{$\displaystyle{p}$(7)\}]
\\\hspace*{8mm}=$\displaystyle{p}$(0)+$\displaystyle{p}$(1)+$\displaystyle{p}$(2)+$\displaystyle{p}$(3)+$\displaystyle{p}$(4)+$\displaystyle{p}$(5)+[$\displaystyle{p}$(6)$\displaystyle{-}$$\displaystyle{p}$(0)$\displaystyle{-}$$\displaystyle{p}$(1)]+[$\displaystyle{p}$(7)$\displaystyle{-}$$\displaystyle{p}$(0)$\displaystyle{-}$$\displaystyle{p}$(1)
\\\hspace*{10.5mm}$\displaystyle{-}$$\displaystyle{p}$(2)$\displaystyle{-}$$\displaystyle{p}$(3)]+[$\displaystyle{p}$(8)$\displaystyle{-}$$\displaystyle{p}$(0)$\displaystyle{-}$$\displaystyle{p}$(1)$\displaystyle{-}$$\displaystyle{p}$(2)$\displaystyle{-}$$\displaystyle{p}$(3)$\displaystyle{-}$$\displaystyle{p}$(4)$\displaystyle{-}$[$\displaystyle{p}$(5)$\displaystyle{-}$\{$\displaystyle{p}$(0)\}
\\\hspace*{10.5mm}$\displaystyle{-}$\{$\displaystyle{p}$(1)\}]]+[$\displaystyle{p}$(9)$\displaystyle{-}$$\displaystyle{p}$(0)$\displaystyle{-}$$\displaystyle{p}$(1)$\displaystyle{-}$$\displaystyle{p}$(2)$\displaystyle{-}$$\displaystyle{p}$(3)$\displaystyle{-}$$\displaystyle{p}$(4)$\displaystyle{-}$[$\displaystyle{p}$(5)$\displaystyle{-}$\{$\displaystyle{p}$(0)\}]$\displaystyle{-}$[$\displaystyle{p}$(6)
\\\hspace*{10.5mm}$\displaystyle{-}$\{$\displaystyle{p}$(0)\}$\displaystyle{-}$\{$\displaystyle{p}$(1)\}$\displaystyle{-}$\{$\displaystyle{p}$(2)\}]$\displaystyle{-}$[$\displaystyle{p}$(7)$\displaystyle{-}$\{$\displaystyle{p}$(0)\}$\displaystyle{-}$\{$\displaystyle{p}$(1)\}$\displaystyle{-}$\{$\displaystyle{p}$(2)\}$\displaystyle{-}$\{$\displaystyle{p}$(3)\}
\\\hspace*{10.5mm}$\displaystyle{-}$\{$\displaystyle{p}$(4)\}]]
\\\hspace*{8mm}=$\displaystyle{p}$(0)+$\displaystyle{p}$(1)+$\displaystyle{p}$(2)+$\displaystyle{p}$(3)+$\displaystyle{p}$(4)+$\displaystyle{p}$(5)+[$\displaystyle{p}$(6)$\displaystyle{-}$$\displaystyle{p}$(0)$\displaystyle{-}$$\displaystyle{p}$(1)]+[$\displaystyle{p}$(7)$\displaystyle{-}$$\displaystyle{p}$(0)$\displaystyle{-}$$\displaystyle{p}$(1)
\\\hspace*{10.5mm}$\displaystyle{-}$$\displaystyle{p}$(2)$\displaystyle{-}$$\displaystyle{p}$(3)]+[$\displaystyle{p}$(8)$\displaystyle{-}$$\displaystyle{p}$(0)$\displaystyle{-}$$\displaystyle{p}$(1)$\displaystyle{-}$$\displaystyle{p}$(2)$\displaystyle{-}$$\displaystyle{p}$(3)$\displaystyle{-}$$\displaystyle{p}$(4)$\displaystyle{-}$[$\displaystyle{p}$(5)$\displaystyle{-}$$\displaystyle{p}$(0)$\displaystyle{-}$$\displaystyle{p}$(1)]]
\\\hspace*{10.5mm}+[$\displaystyle{p}$(9)$\displaystyle{-}$$\displaystyle{p}$(0)$\displaystyle{-}$$\displaystyle{p}$(1)$\displaystyle{-}$$\displaystyle{p}$(2)$\displaystyle{-}$$\displaystyle{p}$(3)$\displaystyle{-}$$\displaystyle{p}$(4)$\displaystyle{-}$[$\displaystyle{p}$(5)$\displaystyle{-}$$\displaystyle{p}$(0)]$\displaystyle{-}$[$\displaystyle{p}$(6)$\displaystyle{-}$$\displaystyle{p}$(0)$\displaystyle{-}$$\displaystyle{p}$(1)
\\\hspace*{10.5mm}$\displaystyle{-}$$\displaystyle{p}$(2)]$\displaystyle{-}$[$\displaystyle{p}$(7)$\displaystyle{-}$$\displaystyle{p}$(0)$\displaystyle{-}$$\displaystyle{p}$(1)$\displaystyle{-}$$\displaystyle{p}$(2)$\displaystyle{-}$$\displaystyle{p}$(3)$\displaystyle{-}$[$\displaystyle{p}$(4)$\displaystyle{-}$\{$\displaystyle{p}$(0)\}]]]
\\\hspace*{8mm}=$\displaystyle{p}$(0)+$\displaystyle{p}$(1)+$\displaystyle{p}$(2)+$\displaystyle{p}$(3)+$\displaystyle{p}$(4)+$\displaystyle{p}$(5)+[$\displaystyle{p}$(6)$\displaystyle{-}$$\displaystyle{p}$(0)$\displaystyle{-}$$\displaystyle{p}$(1)]+[$\displaystyle{p}$(7)$\displaystyle{-}$$\displaystyle{p}$(0)$\displaystyle{-}$$\displaystyle{p}$(1)
\\\hspace*{10.5mm}$\displaystyle{-}$$\displaystyle{p}$(2)$\displaystyle{-}$$\displaystyle{p}$(3)]+[$\displaystyle{p}$(8)$\displaystyle{-}$$\displaystyle{p}$(0)$\displaystyle{-}$$\displaystyle{p}$(1)$\displaystyle{-}$$\displaystyle{p}$(2)$\displaystyle{-}$$\displaystyle{p}$(3)$\displaystyle{-}$$\displaystyle{p}$(4)$\displaystyle{-}$[$\displaystyle{p}$(5)$\displaystyle{-}$$\displaystyle{p}$(0)$\displaystyle{-}$$\displaystyle{p}$(1)]]
\\\hspace*{10.5mm}+[$\displaystyle{p}$(9)$\displaystyle{-}$$\displaystyle{p}$(0)$\displaystyle{-}$$\displaystyle{p}$(1)$\displaystyle{-}$$\displaystyle{p}$(2)$\displaystyle{-}$$\displaystyle{p}$(3)$\displaystyle{-}$$\displaystyle{p}$(4)$\displaystyle{-}$[$\displaystyle{p}$(5)$\displaystyle{-}$$\displaystyle{p}$(0)]$\displaystyle{-}$[$\displaystyle{p}$(6)$\displaystyle{-}$$\displaystyle{p}$(0)$\displaystyle{-}$$\displaystyle{p}$(1)
\\\hspace*{10.5mm}$\displaystyle{-}$$\displaystyle{p}$(2)]$\displaystyle{-}$[$\displaystyle{p}$(7)$\displaystyle{-}$$\displaystyle{p}$(0)$\displaystyle{-}$$\displaystyle{p}$(1)$\displaystyle{-}$$\displaystyle{p}$(2)$\displaystyle{-}$$\displaystyle{p}$(3)$\displaystyle{-}$[$\displaystyle{p}$(4)$\displaystyle{-}$$\displaystyle{p}$(0)]]]
\\\hspace*{8mm}=42.\hfill(3)

One can note that the equation above has the feature: self-similarity, of course, this is the property of the fractal.

On the one hand $\displaystyle{p(n)}$ can be ``zoomed in'' on a form as $\displaystyle{p(n)}$\,=\,$\displaystyle{\sum_{i=0}^{n-1}}$$\displaystyle{\{p(i)\}}$, on the other hand, in fact, for $\displaystyle{p(n)}$, the last term \{$\displaystyle{p(n-1)}$\} in $\displaystyle{p(n)}$ can be mapped as the sum of $\displaystyle{n}$ integer 1, there is only one way to express \{$\displaystyle{p(n-1)}$\}, so the last term can be 1. Thus $\displaystyle{p(n)}$ can be ``zoomed in'' on a form as $\displaystyle{p(n)}$\,=\,$\displaystyle{\sum_{i=0}^{n-2}}$$\displaystyle{\{p(i)\}}$\,+\,1.

Also, for $\displaystyle{p}$(10), it can be ``zoomed in'' entirely with 4 steps, we have\vspace{1ex}
\\$\displaystyle{p}$(10)=$\displaystyle{\sum_{i=0}^{8}}$$\displaystyle{\{p(i)\}}$+1
\\\hspace*{8mm}=$\displaystyle{p}$(0)+$\displaystyle{p}$(1)+$\displaystyle{p}$(2)+$\displaystyle{p}$(3)+$\displaystyle{p}$(4)+$\displaystyle{p}$(5)+[$\displaystyle{p}$(6)$\displaystyle{-}$\{$\displaystyle{p}$(0)\}$\displaystyle{-}$\{$\displaystyle{p}$(1)\}]+[$\displaystyle{p}$(7)$\displaystyle{-}$\{$\displaystyle{p}$(0)\}
\\\hspace*{10.5mm}$\displaystyle{-}$\{$\displaystyle{p}$(1)\}$\displaystyle{-}$\{$\displaystyle{p}$(2)\}$\displaystyle{-}$\{$\displaystyle{p}$(3)\}]+[$\displaystyle{p}$(8)$\displaystyle{-}$\{$\displaystyle{p}$(0)\}$\displaystyle{-}$\{$\displaystyle{p}$(1)\}$\displaystyle{-}$\{$\displaystyle{p}$(2)\}$\displaystyle{-}$\{$\displaystyle{p}$(3)\}
\\\hspace*{10.5mm}$\displaystyle{-}$\{$\displaystyle{p}$(4)\}$\displaystyle{-}$\{$\displaystyle{p}$(5)\}]+1
\\\hspace*{8mm}=$\displaystyle{p}$(0)+$\displaystyle{p}$(1)+$\displaystyle{p}$(2)+$\displaystyle{p}$(3)+$\displaystyle{p}$(4)+$\displaystyle{p}$(5)+[$\displaystyle{p}$(6)$\displaystyle{-}$$\displaystyle{p}$(0)$\displaystyle{-}$$\displaystyle{p}$(1)]+[$\displaystyle{p}$(7)$\displaystyle{-}$$\displaystyle{p}$(0)$\displaystyle{-}$$\displaystyle{p}$(1)
\\\hspace*{10.5mm}$\displaystyle{-}$$\displaystyle{p}$(2)$\displaystyle{-}$$\displaystyle{p}$(3)]+[$\displaystyle{p}$(8)$\displaystyle{-}$$\displaystyle{p}$(0)$\displaystyle{-}$$\displaystyle{p}$(1)$\displaystyle{-}$$\displaystyle{p}$(2)$\displaystyle{-}$$\displaystyle{p}$(3)$\displaystyle{-}$$\displaystyle{p}$(4)$\displaystyle{-}$[$\displaystyle{p}$(5)$\displaystyle{-}$\{$\displaystyle{p}$(0)\}
\\\hspace*{10.5mm}$\displaystyle{-}$\{$\displaystyle{p}$(1)\}]]+1
\\\hspace*{8mm}=$\displaystyle{p}$(0)+$\displaystyle{p}$(1)+$\displaystyle{p}$(2)+$\displaystyle{p}$(3)+$\displaystyle{p}$(4)+$\displaystyle{p}$(5)+[$\displaystyle{p}$(6)$\displaystyle{-}$$\displaystyle{p}$(0)$\displaystyle{-}$$\displaystyle{p}$(1)]+[$\displaystyle{p}$(7)$\displaystyle{-}$$\displaystyle{p}$(0)$\displaystyle{-}$$\displaystyle{p}$(1)
\\\hspace*{10.5mm}$\displaystyle{-}$$\displaystyle{p}$(2)$\displaystyle{-}$$\displaystyle{p}$(3)]+[$\displaystyle{p}$(8)$\displaystyle{-}$$\displaystyle{p}$(0)$\displaystyle{-}$$\displaystyle{p}$(1)$\displaystyle{-}$$\displaystyle{p}$(2)$\displaystyle{-}$$\displaystyle{p}$(3)$\displaystyle{-}$$\displaystyle{p}$(4)$\displaystyle{-}$[$\displaystyle{p}$(5)$\displaystyle{-}$$\displaystyle{p}$(0)$\displaystyle{-}$$\displaystyle{p}$(1)]]+1
\\\hspace*{8mm}=42.\hfill(4)

It is easy to find that $\displaystyle{p(n)}$ also can be ``zoomed in'' on a form as \\$\displaystyle{p(n)}$\,=\,$\displaystyle{\sum_{i=0}^{n-3}}$$\displaystyle{\{p(i)\}}$\,+\,$\lfloor n/2 \rfloor$\,+\,1 if we have gotten the content which will be discussed later. Also $\displaystyle{p}$(10) can be ``zoomed in'' entirely with 3 steps, we have\vspace{1ex}
\\$\displaystyle{p}$(10)=$\displaystyle{\sum_{i=0}^{7}}$$\displaystyle{\{p(i)\}}$+$\lfloor 10/2 \rfloor$+1
\\\hspace*{8mm}=$\displaystyle{p}$(0)+$\displaystyle{p}$(1)+$\displaystyle{p}$(2)+$\displaystyle{p}$(3)+$\displaystyle{p}$(4)+$\displaystyle{p}$(5)+[$\displaystyle{p}$(6)$\displaystyle{-}$\{$\displaystyle{p}$(0)\}$\displaystyle{-}$\{$\displaystyle{p}$(1)\}]+[$\displaystyle{p}$(7)$\displaystyle{-}$\{$\displaystyle{p}$(0)\}
\\\hspace*{10.5mm}$\displaystyle{-}$\{$\displaystyle{p}$(1)\}$\displaystyle{-}$\{$\displaystyle{p}$(2)\}$\displaystyle{-}$\{$\displaystyle{p}$(3)\}]+10/2+1
\\\hspace*{8mm}=$\displaystyle{p}$(0)+$\displaystyle{p}$(1)+$\displaystyle{p}$(2)+$\displaystyle{p}$(3)+$\displaystyle{p}$(4)+$\displaystyle{p}$(5)+[$\displaystyle{p}$(6)$\displaystyle{-}$$\displaystyle{p}$(0)$\displaystyle{-}$$\displaystyle{p}$(1)]+[$\displaystyle{p}$(7)$\displaystyle{-}$$\displaystyle{p}$(0)$\displaystyle{-}$$\displaystyle{p}$(1)
\\\hspace*{10.5mm}$\displaystyle{-}$$\displaystyle{p}$(2)$\displaystyle{-}$$\displaystyle{p}$(3)]+10/2+1
\\\hspace*{8mm}=42.\hfill(5)

In the fact that $\displaystyle{p(n)}$ can be ``zoomed in'' entirely with fewer steps if we can find some proper functions and exchange the terms in $\displaystyle{p(n)}$ for these functions.

\section{The application of the generator}

\noindent Now let us reach an agreement that the the main theorem in the section above is called the generator. As an example, let us use the generator to calculate $\displaystyle{p(n)}$ in a finite range. Here the last term \{$\displaystyle{p(n-1)}$\} can be mapped as the sum of $\displaystyle{(n-1)}$ integer 1 to add up, there is only one way to express {\{$\displaystyle{p(n-1)}$\}}, so the last term is 1.

Let us assume 2 $\leq$ $\displaystyle{n}$ $\leq$ 12, at first we have \vspace{0.5ex}

$\displaystyle{p}$($\displaystyle{n}$)=\{$\displaystyle{p}$($\displaystyle{n}$$\displaystyle{-}$$\displaystyle{n}$)\}+ $\cdots$ +\{$\displaystyle{p}$($\lfloor n/2 \rfloor$)\}+\{$\displaystyle{p}$($\displaystyle{n}$$\displaystyle{-}$$\displaystyle{x}$)\}+ $\cdots$ +\{$\displaystyle{p}$($\displaystyle{n}$$\displaystyle{-}$2)\}+1.\vspace{0.5ex}\hfill(6)

We should calculate $\displaystyle{x}$ first. \{$\displaystyle{p(n-x)}$\} is the first parcel that will generate its child-parcel(s), $\displaystyle{n-x}$ $>$ $\displaystyle{n/2}$, then $\displaystyle{x}$ $<$ $\displaystyle{n}$ $\displaystyle{-}$ $\displaystyle{n/2}$, therefore, $\displaystyle{x}$ $<$ 6.

The tabs in the equation from left to right are arranged from small to large, therefore $\displaystyle{x}$ should has the largest value in the range. Thus we know $\displaystyle{x}$ = 5.

Using the generator, we know the quantity of child-parcel(s) that a parcel will generate, also we know the first child-parcel should be $\displaystyle{p(0)}$ and the child-parcel(s) are arranged by order according to their tabs. The quantity of the child-parcel(s) that a parcel will generate is given by the generator, therefore the tab of the last child-parcel equals the quantity of the child-parcel(s) minus 1. $\displaystyle{p}$($\lfloor n/2 \rfloor$) = $\displaystyle{p(n-6)}$ because the number $\lfloor n/2 \rfloor$ in $\displaystyle{p}$($\lfloor n/2 \rfloor$) and the number $\displaystyle{(n-5)}$ in $\displaystyle{p(n-5)}$ should be continuous, we have \vspace{1ex}
\\$\displaystyle{p}$($\displaystyle{n}$)=\{$\displaystyle{p}$($\displaystyle{n}$$\displaystyle{-}$$\displaystyle{n}$)\}+ $\cdots$ +\{$\displaystyle{p}$($\lfloor n/2 \rfloor$)\}+\{$\displaystyle{p}$($\displaystyle{n}$$\displaystyle{-}$5)\}+ $\cdots$ +\{$\displaystyle{p}$($\displaystyle{n}$$\displaystyle{-}$2)\}+1\vspace{1ex}
\\=$\displaystyle{p}$($\displaystyle{n}$$\displaystyle{-}$$\displaystyle{n}$)+ $\cdots$ +$\displaystyle{p}$($\displaystyle{n}$$\displaystyle{-}$6)
\\+[$\displaystyle{p}$($\displaystyle{n}$$\displaystyle{-}$5)$\displaystyle{-}$$\displaystyle{p}$($\displaystyle{n}$$\displaystyle{-}$12)$\displaystyle{-}$$\displaystyle{p}$($\displaystyle{n}$$\displaystyle{-}$11)]
\\+[$\displaystyle{p}$($\displaystyle{n}$$\displaystyle{-}$4)$\displaystyle{-}$$\displaystyle{p}$($\displaystyle{n}$$\displaystyle{-}$12)$\displaystyle{-}$$\displaystyle{p}$($\displaystyle{n}$$\displaystyle{-}$11)$\displaystyle{-}$$\displaystyle{p}$($\displaystyle{n}$$\displaystyle{-}$10)$\displaystyle{-}$$\displaystyle{p}$($\displaystyle{n}$$\displaystyle{-}$9)]
\\+[$\displaystyle{p}$($\displaystyle{n}$$\displaystyle{-}$3)$\displaystyle{-}$$\displaystyle{p}$($\displaystyle{n}$$\displaystyle{-}$12)$\displaystyle{-}$$\displaystyle{p}$($\displaystyle{n}$$\displaystyle{-}$11)$\displaystyle{-}$$\displaystyle{p}$($\displaystyle{n}$$\displaystyle{-}$10)$\displaystyle{-}$$\displaystyle{p}$($\displaystyle{n}$$\displaystyle{-}$9)$\displaystyle{-}$$\displaystyle{p}$($\displaystyle{n}$$\displaystyle{-}$8)$\displaystyle{-}$[$\displaystyle{p}$($\displaystyle{n}$$\displaystyle{-}$7)
\\$\displaystyle{-}$$\displaystyle{p}$($\displaystyle{n}$$\displaystyle{-}$12)]]
\\+[$\displaystyle{p}$($\displaystyle{n}$$\displaystyle{-}$2)$\displaystyle{-}$$\displaystyle{p}$($\displaystyle{n}$$\displaystyle{-}$12)$\displaystyle{-}$$\displaystyle{p}$($\displaystyle{n}$$\displaystyle{-}$11)$\displaystyle{-}$$\displaystyle{p}$($\displaystyle{n}$$\displaystyle{-}$10)$\displaystyle{-}$$\displaystyle{p}$($\displaystyle{n}$$\displaystyle{-}$9)$\displaystyle{-}$$\displaystyle{p}$($\displaystyle{n}$$\displaystyle{-}$8)$\displaystyle{-}$$\displaystyle{p}$($\displaystyle{n}$$\displaystyle{-}$7) \\$\displaystyle{-}$[$\displaystyle{p}$($\displaystyle{n}$$\displaystyle{-}$6)$\displaystyle{-}$$\displaystyle{p}$($\displaystyle{n}$$\displaystyle{-}$12)$\displaystyle{-}$$\displaystyle{p}$($\displaystyle{n}$$\displaystyle{-}$11)]$\displaystyle{-}$[$\displaystyle{p}$($\displaystyle{n}$$\displaystyle{-}$5)$\displaystyle{-}$$\displaystyle{p}$($\displaystyle{n}$$\displaystyle{-}$12)$\displaystyle{-}$$\displaystyle{p}$($\displaystyle{n}$$\displaystyle{-}$11)
\\$\displaystyle{-}$$\displaystyle{p}$($\displaystyle{n}$$\displaystyle{-}$10)$\displaystyle{-}$$\displaystyle{p}$($\displaystyle{n}$$\displaystyle{-}$9)]]+1.
\vspace{1ex}\hfill(7)

Given a nonnegative integer $\displaystyle{n}$, now let us use the generator to calculate $\displaystyle{p(n)}$ directly.

\noindent  \emph{Example:} Calculate $\displaystyle{p(11)}$ and $\displaystyle{p(12)}$ by using the generator.

\vspace{1ex} \noindent $\displaystyle{p}$(11)=$\displaystyle{p}$(0)+$\displaystyle{p}$(1)+$\displaystyle{p}$(2)+$\displaystyle{p}$(3)+$\displaystyle{p}$(4)+$\displaystyle{p}$(5)+[$\displaystyle{p}$(6)$\displaystyle{-}$$\displaystyle{p}$(0)]+[$\displaystyle{p}$(7)$\displaystyle{-}$$\displaystyle{p}$(0)$\displaystyle{-}$$\displaystyle{p}$(1)$\displaystyle{-}$$\displaystyle{p}$(2)]
\\
\hspace*{10.5mm}+[$\displaystyle{p}$(8)$\displaystyle{-}$$\displaystyle{p}$(0)$\displaystyle{-}$$\displaystyle{p}$(1)$\displaystyle{-}$$\displaystyle{p}$(2)$\displaystyle{-}$$\displaystyle{p}$(3)$\displaystyle{-}$$\displaystyle{p}$(4)]+[$\displaystyle{p}$(9)$\displaystyle{-}$$\displaystyle{p}$(0)$\displaystyle{-}$$\displaystyle{p}$(1)$\displaystyle{-}$$\displaystyle{p}$(2)$\displaystyle{-}$$\displaystyle{p}$(3)
\\
\hspace*{10.5mm}$\displaystyle{-}$$\displaystyle{p}$(4)$\displaystyle{-}$[$\displaystyle{p}$(5)$\displaystyle{-}$$\displaystyle{p}$(0)]$\displaystyle{-}$[$\displaystyle{p}$(6)$\displaystyle{-}$$\displaystyle{p}$(0)$\displaystyle{-}$$\displaystyle{p}$(1)$\displaystyle{-}$$\displaystyle{p}$(2)]]+1
\\
\hspace*{8mm}=56,\vspace{1ex}\hfill(8)
\noindent
\\$\displaystyle{p}$(12)=$\displaystyle{p}$(0)+$\displaystyle{p}$(1)+$\displaystyle{p}$(2)+$\displaystyle{p}$(3)+$\displaystyle{p}$(4)+$\displaystyle{p}$(5)+$\displaystyle{p}$(6)+[$\displaystyle{p}$(7)$\displaystyle{-}$$\displaystyle{p}$(0)$\displaystyle{-}$$\displaystyle{p}$(1)]+[$\displaystyle{p}$(8)$\displaystyle{-}$$\displaystyle{p}$(0)
\\
\hspace*{10.5mm}$\displaystyle{-}$$\displaystyle{p}$(1)$\displaystyle{-}$$\displaystyle{p}$(2)$\displaystyle{-}$$\displaystyle{p}$(3)]+[$\displaystyle{p}$(9)$\displaystyle{-}$$\displaystyle{p}$(0)$\displaystyle{-}$$\displaystyle{p}$(1)$\displaystyle{-}$$\displaystyle{p}$(2)$\displaystyle{-}$$\displaystyle{p}$(3)$\displaystyle{-}$$\displaystyle{p}$(4)$\displaystyle{-}$[$\displaystyle{p}$(5)$\displaystyle{-}$$\displaystyle{p}$(0)]]
\\
\hspace*{10.5mm}+[$\displaystyle{p}$(10)$\displaystyle{-}$$\displaystyle{p}$(0)$\displaystyle{-}$$\displaystyle{p}$(1)$\displaystyle{-}$$\displaystyle{p}$(2)$\displaystyle{-}$$\displaystyle{p}$(3)$\displaystyle{-}$$\displaystyle{p}$(4)$\displaystyle{-}$$\displaystyle{p}$(5)$\displaystyle{-}$[$\displaystyle{p}$(6)$\displaystyle{-}$$\displaystyle{p}$(0)$\displaystyle{-}$$\displaystyle{p}$(1)]$\displaystyle{-}$[$\displaystyle{p}$(7)
\\
\hspace*{10.5mm}$\displaystyle{-}$$\displaystyle{p}$(0)$\displaystyle{-}$$\displaystyle{p}$(1)$\displaystyle{-}$$\displaystyle{p}$(2)$\displaystyle{-}$$\displaystyle{p}$(3)]]+1
\\
\hspace*{8mm}=77.\vspace{1ex}\hfill(9)

Let us calculate $\displaystyle{p(n-1)}$ when 2 $\leq$ $\displaystyle{n}$ $\leq$ 12, this is a preparing work for the next chapter. The process of analysis is as same as $\displaystyle{p(n)}$ above. \vspace{1ex}
\\$\displaystyle{p}$($\displaystyle{n}$$\displaystyle{-}$1)=\{$\displaystyle{p}$(($\displaystyle{n}$$\displaystyle{-}$1)$\displaystyle{-}$($\displaystyle{n}$$\displaystyle{-}$1))\}+ $\cdots$ +\{$\displaystyle{p}$($\lfloor (n-1)/2 \rfloor$)\}+\{$\displaystyle{p}$($\displaystyle{n}$$\displaystyle{-}$1$\displaystyle{-}$5)\}+ $\cdots$ +
\\\{$\displaystyle{p}$($\displaystyle{n}$$\displaystyle{-}$1$\displaystyle{-}$2)\}+1\vspace{1ex}
\\=$\displaystyle{p}$($\displaystyle{n}$$\displaystyle{-}$$\displaystyle{n}$)+ $\cdots$ +$\displaystyle{p}$($\displaystyle{n}$$\displaystyle{-}$7)
\\+[$\displaystyle{p}$($\displaystyle{n}$$\displaystyle{-}$6)$\displaystyle{-}$$\displaystyle{p}$($\displaystyle{n}$$\displaystyle{-}$12)]
\\+[$\displaystyle{p}$($\displaystyle{n}$$\displaystyle{-}$5)$\displaystyle{-}$$\displaystyle{p}$($\displaystyle{n}$$\displaystyle{-}$12)$\displaystyle{-}$$\displaystyle{p}$($\displaystyle{n}$$\displaystyle{-}$11)$\displaystyle{-}$$\displaystyle{p}$($\displaystyle{n}$$\displaystyle{-}$10)] \\+[$\displaystyle{p}$($\displaystyle{n}$$\displaystyle{-}$4)$\displaystyle{-}$$\displaystyle{p}$($\displaystyle{n}$$\displaystyle{-}$12)$\displaystyle{-}$$\displaystyle{p}$($\displaystyle{n}$$\displaystyle{-}$11)$\displaystyle{-}$$\displaystyle{p}$($\displaystyle{n}$$\displaystyle{-}$10)$\displaystyle{-}$$\displaystyle{p}$($\displaystyle{n}$$\displaystyle{-}$9)$\displaystyle{-}$$\displaystyle{p}$($\displaystyle{n}$$\displaystyle{-}$8)]
\\+[$\displaystyle{p}$($\displaystyle{n}$$\displaystyle{-}$3)$\displaystyle{-}$$\displaystyle{p}$($\displaystyle{n}$$\displaystyle{-}$12)$\displaystyle{-}$$\displaystyle{p}$($\displaystyle{n}$$\displaystyle{-}$11)$\displaystyle{-}$$\displaystyle{p}$($\displaystyle{n}$$\displaystyle{-}$10)$\displaystyle{-}$$\displaystyle{p}$($\displaystyle{n}$$\displaystyle{-}$9)$\displaystyle{-}$$\displaystyle{p}$($\displaystyle{n}$$\displaystyle{-}$8)$\displaystyle{-}$[$\displaystyle{p}$($\displaystyle{n}$$\displaystyle{-}$7)
\\$\displaystyle{-}$$\displaystyle{p}$($\displaystyle{n}$$\displaystyle{-}$12)] \\$\displaystyle{-}$[$\displaystyle{p}$($\displaystyle{n}$$\displaystyle{-}$6)$\displaystyle{-}$$\displaystyle{p}$($\displaystyle{n}$$\displaystyle{-}$12)$\displaystyle{-}$$\displaystyle{p}$($\displaystyle{n}$$\displaystyle{-}$11)$\displaystyle{-}$$\displaystyle{p}$($\displaystyle{n}$$\displaystyle{-}$10)]]+1.
\hfill(10)

\section{Deducing some recurrence equations}

\noindent Now let us consider the equation (7) and (10), let $\displaystyle{p(n)-p(n-1)}$, we have \vspace{1ex}
\\$\displaystyle{p(n)=p(n-1)+p(n-2)+p(n-12)-p(n-5)-p(n-7)}$.                                      \vspace{1ex}\hfill(11)

That is the pentagonal number theorem \textsuperscript{\cite{3}-\cite{8}} when 2 $\leq$ $\displaystyle{n}$ $\leq$ 12.

Let us change the tail of the equation (7) or (10) to get some new results.

It is easy to prove that when $\displaystyle{n}$ is an even number, the count of the integer partitions which 2 is the largest number equals $\displaystyle{n/2}$, when $\displaystyle{n}$ is an odd number, it equals $\displaystyle{(n-1)/2}$.

$\displaystyle\{p(n-2)\}$ can be exchanged for $\lfloor n/2 \rfloor$ because \{$\displaystyle{p(n-2)}$\} equals the count of the integer partitions which 2 is the largest number.

Therefore, $\displaystyle{p}$(11) and $\displaystyle{p}$(12) can be calculated as follows:\vspace{1ex}

\noindent
$\displaystyle{p}$(11)=$\displaystyle{p}$(0)+$\displaystyle{p}$(1)+$\displaystyle{p}$(2)+$\displaystyle{p}$(3)+$\displaystyle{p}$(4)+$\displaystyle{p}$(5)+[$\displaystyle{p}$(6)$\displaystyle{-}$$\displaystyle{p}$(0)]+[$\displaystyle{p}$(7)$\displaystyle{-}$$\displaystyle{p}$(0)$\displaystyle{-}$$\displaystyle{p}$(1)$\displaystyle{-}$$\displaystyle{p}$(2)]
\\
\hspace*{10.5mm}+[$\displaystyle{p}$(8)$\displaystyle{-}$$\displaystyle{p}$(0)$\displaystyle{-}$$\displaystyle{p}$(1)$\displaystyle{-}$$\displaystyle{p}$(2)$\displaystyle{-}$$\displaystyle{p}$(3)$\displaystyle{-}$$\displaystyle{p}$(4)]+(11$\displaystyle{-}$1)/2+1
\\
\hspace*{8mm}=56,
\vspace{2ex}\hfill(12)
\noindent
\\$\displaystyle{p}$(12)=$\displaystyle{p}$(0)+$\displaystyle{p}$(1)+$\displaystyle{p}$(2)+$\displaystyle{p}$(3)+$\displaystyle{p}$(4)+$\displaystyle{p}$(5)+$\displaystyle{p}$(6)+[$\displaystyle{p}$(7)$\displaystyle{-}$$\displaystyle{p}$(0)$\displaystyle{-}$$\displaystyle{p}$(1)]+[$\displaystyle{p}$(8)$\displaystyle{-}$$\displaystyle{p}$(0)
\\
\hspace*{10.5mm}$\displaystyle{-}$$\displaystyle{p}$(1)$\displaystyle{-}$$\displaystyle{p}$(2)$\displaystyle{-}$$\displaystyle{p}$(3)]+[$\displaystyle{p}$(9)$\displaystyle{-}$$\displaystyle{p}$(0)$\displaystyle{-}$$\displaystyle{p}$(1)$\displaystyle{-}$$\displaystyle{p}$(2)$\displaystyle{-}$$\displaystyle{p}$(3)$\displaystyle{-}$$\displaystyle{p}$(4)$\displaystyle{-}$[$\displaystyle{p}$(5)$\displaystyle{-}$$\displaystyle{p}$(0)]]+
\\
\hspace*{10.5mm}12/2+1
\\
\hspace*{8mm}=77.\vspace{3ex}\hfill(13)

Now we can see that if we let $\displaystyle{p(n)}$\,=\,$\displaystyle{\sum_{i=0}^{n-3}}$$\displaystyle{\{p(i)\}}$\,+\,$\lfloor n/2 \rfloor$\,+\,1 and let \\$\displaystyle{p(n-1)}$\,=\,$\displaystyle{\sum_{i=0}^{n-3}}$$\displaystyle{\{p(i)\}}$\,+\,1, in fact, in a larger range, when 2 $\leq$ $\displaystyle{n}$ $\leq$ 24, after being ``zoomed in'' entirely both for $\displaystyle{p(n)}$ and $\displaystyle{p(n-1)}$, making $\displaystyle{p(n)-p(n-1)}$, we have\vspace{3.5ex}
\\$\displaystyle{p(n)=p(n-1)+p(n-6)+p(n-8)+p(n-20)+p(n-21)+p(n-22)}$
\\$\displaystyle{+2p(n-23)+2p(n-24)-p(n-11)-2p(n-13)-p(n-14)-p(n-15)}$
\\$\displaystyle{-p(n-16)-p(n-17)+(n-k)/2}$
\hfill($\displaystyle{k}$ = 0 if $\displaystyle{n}$ is even, $\displaystyle{k}$ = 1 if $\displaystyle{n}$ is odd).\vspace{3.5ex}
\hfill(14)

Let
$\displaystyle{p(n)}$\,=\,$\displaystyle{\sum_{i=0}^{n-3}}$$\displaystyle{\{p(i)\}}$\,+\,$\lfloor n/2 \rfloor$\,+\,1 and let \\ $\displaystyle{p(n-1)}$\,=\,$\displaystyle{\sum_{i=0}^{n-4}}$$\displaystyle{\{p(i)\}}$\,+\,$\lfloor (n-1)/2 \rfloor$\,+\,1, when 2 $\leq$ $\displaystyle{n}$ $\leq$ 24, after being ``zoomed in'' entirely both for $\displaystyle{p(n)}$ and $\displaystyle{p(n-1)}$, making $\displaystyle{p(n)-p(n-1)}$, we have\vspace{3.5ex}
\\$\displaystyle{p(n)=p(n-1)+p(n-3)+p(n-12)+p(n-14)+p(n-16)+p(n-18)}$
\\$\displaystyle{+p(n-20)-p(n-7)-p(n-9)-p(n-11)-p(n-13)+k}$
\hfill($\displaystyle{k}$ = 1 if $\displaystyle{n}$ is even, $\displaystyle{k}$ = 0 if $\displaystyle{n}$ is odd).\vspace{3.5ex}
\hfill(15)

Let $\displaystyle{p(n)}$\,=\,$\displaystyle{\sum_{i=0}^{n-1}}$$\displaystyle{\{p(i)\}}$ and let $\displaystyle{p(n-1)}$\,=\,$\displaystyle{\sum_{i=0}^{n-2}}$$\displaystyle{\{p(i)\}}$, when 2 $\leq$ $\displaystyle{n}$ $\leq$ 24, after being ``zoomed in'' entirely both for $\displaystyle{p(n)}$ and $\displaystyle{p(n-1)}$, making $\displaystyle{p(n)-p(n-1)}$, we have\vspace{3.5ex}
\\$\displaystyle{p(n)=2p(n-1)+p(n-6)+p(n-8)+p(n-12)+p(n-15)-p(n-3)}$
\\$\displaystyle{-p(n-5)-p(n-7)-p(n-13)-p(n-16)-p(n-24)}$.\vspace{3.5ex}
\hfill(16)

Let $\displaystyle{p(n)}$\,=\,$\displaystyle{\sum_{i=0}^{n-2}}$$\displaystyle{\{p(i)\}}$\,+\,1 and let $\displaystyle{p(n-1)}$\,=\,$\displaystyle{\sum_{i=0}^{n-3}}$$\displaystyle{\{p(i)\}}$\,+\,1, also when 2 $\leq$ $\displaystyle{n}$ $\leq$ 24, after being ``zoomed in'' entirely both for $\displaystyle{p(n)}$ and $\displaystyle{p(n-1)}$, making $\displaystyle{p(n)-p(n-1)}$, we have\vspace{3.5ex}
\\$\displaystyle{p(n)=p(n-1)+p(n-2)+p(n-12)+p(n-15)-p(n-5)-p(n-7)}$
\\$\displaystyle{-p(n-22)}$  \,\,\,\,\,\,\,\, (the pentagonal number theorem).\vspace{3.5ex}
\hfill(17)

\section{Conclusions}
\noindent All the recurrence equations appearing above are deduced by the generator which has been deduced in this paper, the pentagonal number theorem is a special case between them, in fact, we can get more recurrence equations.

There is an inference that the count of the kinds of the recurrence equations like above which form is $\displaystyle{F(p(n),p(n-1), \cdots ,p(0))+G(n)=0}$ should be infinite, where $\displaystyle{F(x_0,x_1, \cdots ,x_n)}$ is a simple equation and $\displaystyle{G(n)}$ is a function of $\displaystyle{n}$. It is because of that the count of the terms that the equation $\displaystyle{p(n)}$ has can be infinite if $\displaystyle{n}$ is infinite, in the first step of ``zooming in'' $\displaystyle{p(n)}$, every term of the second half of the equation $\displaystyle{p(n)}$ can be exchanged for a function $\displaystyle{G(n)}$ freely.

Also, formula (14) implies every natural number $\displaystyle{n}$ can be expressed by the sum of positive and negative $\displaystyle{p(i)}$, where 0 $\leq$ $\displaystyle{i}$ $\leq$ $\displaystyle{n}$.

In fact, when 1 $\leq$ $\displaystyle{n}$ $\leq$ 24, we have\vspace{3.5ex}
\\$\displaystyle{n=2[p(n)-p(n-1)-p(n-6)-p(n-8)-p(n-20)-p(n-21)-p(n-22)}$
\\$\displaystyle{-2p(n-23)-2p(n-24)+p(n-11)+2p(n-13)+p(n-14)+p(n-15)}$
\\$\displaystyle{+p(n-16)+p(n-17)]+k}$
\hfill($\displaystyle{k}$ = 0 if $\displaystyle{n}$ is even, $\displaystyle{k}$ = 1 if $\displaystyle{n}$ is odd).\vspace{3.5ex}
\hfill(18)

The formula above may lead us to realise the area unknown.


\bigskip

\noindent \textsc{School of Computer and Information Engineering, Tianjin Chengjian University(The Original Name: Tianjin Institute of Urban Construction), Tianjin 300384, P.R. China}\\
The Author's Chinese Identity Card Number: \textsf{120103198410084517} \\
\textit{E-mail address}: hsmengzh@hotmail.com
\\
\\
The photographs of the author's Chinese identity card, the author's Chinese social security card and the author's individual photos are also shown in the following pages.

\centerline{\includegraphics[width=6cm]{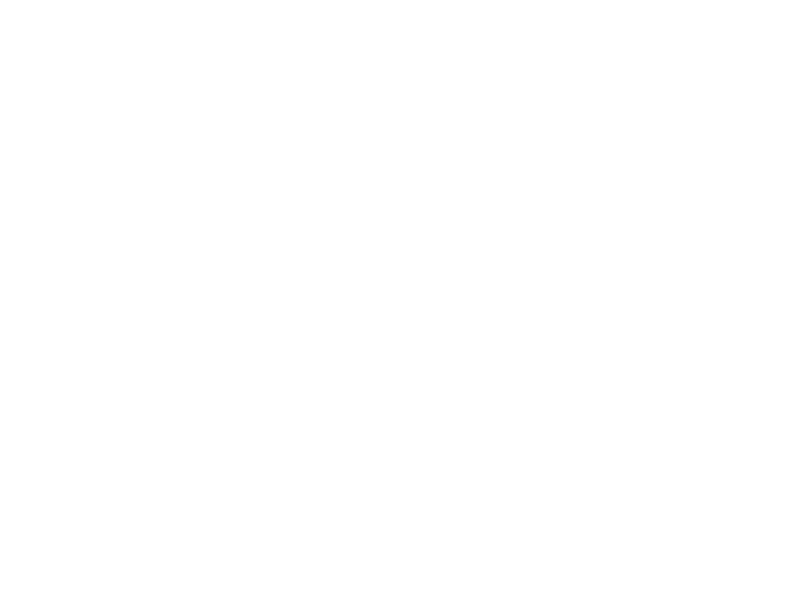}}

\begin{center}
   \includegraphics[width=12cm]{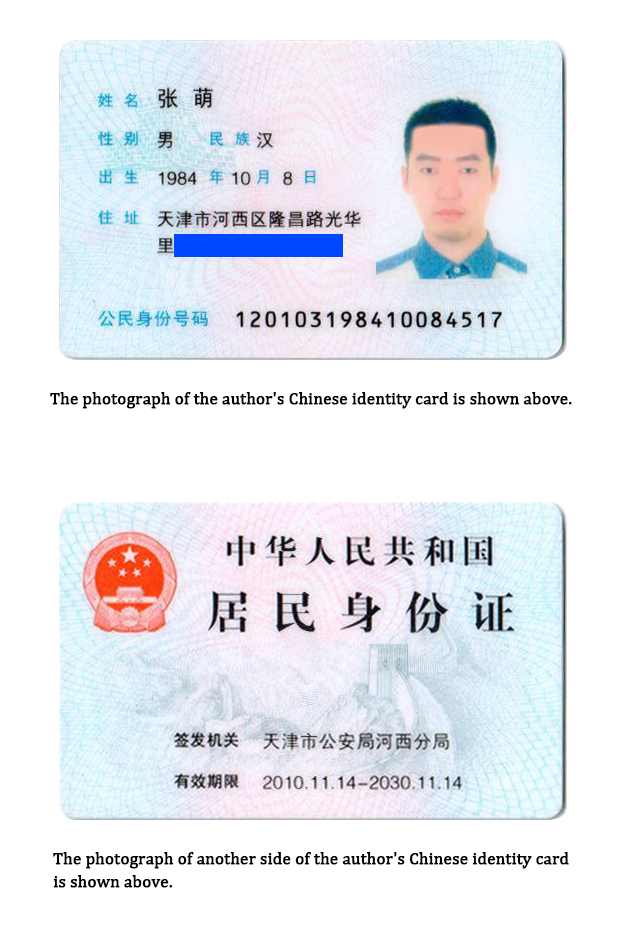}
   \includegraphics[width=12cm]{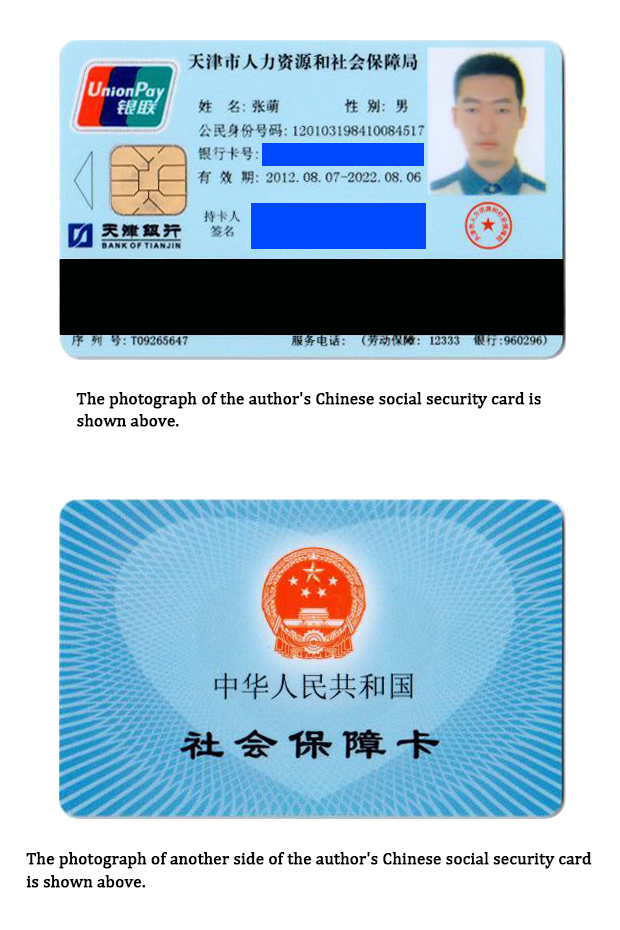}
   \includegraphics[width=12cm]{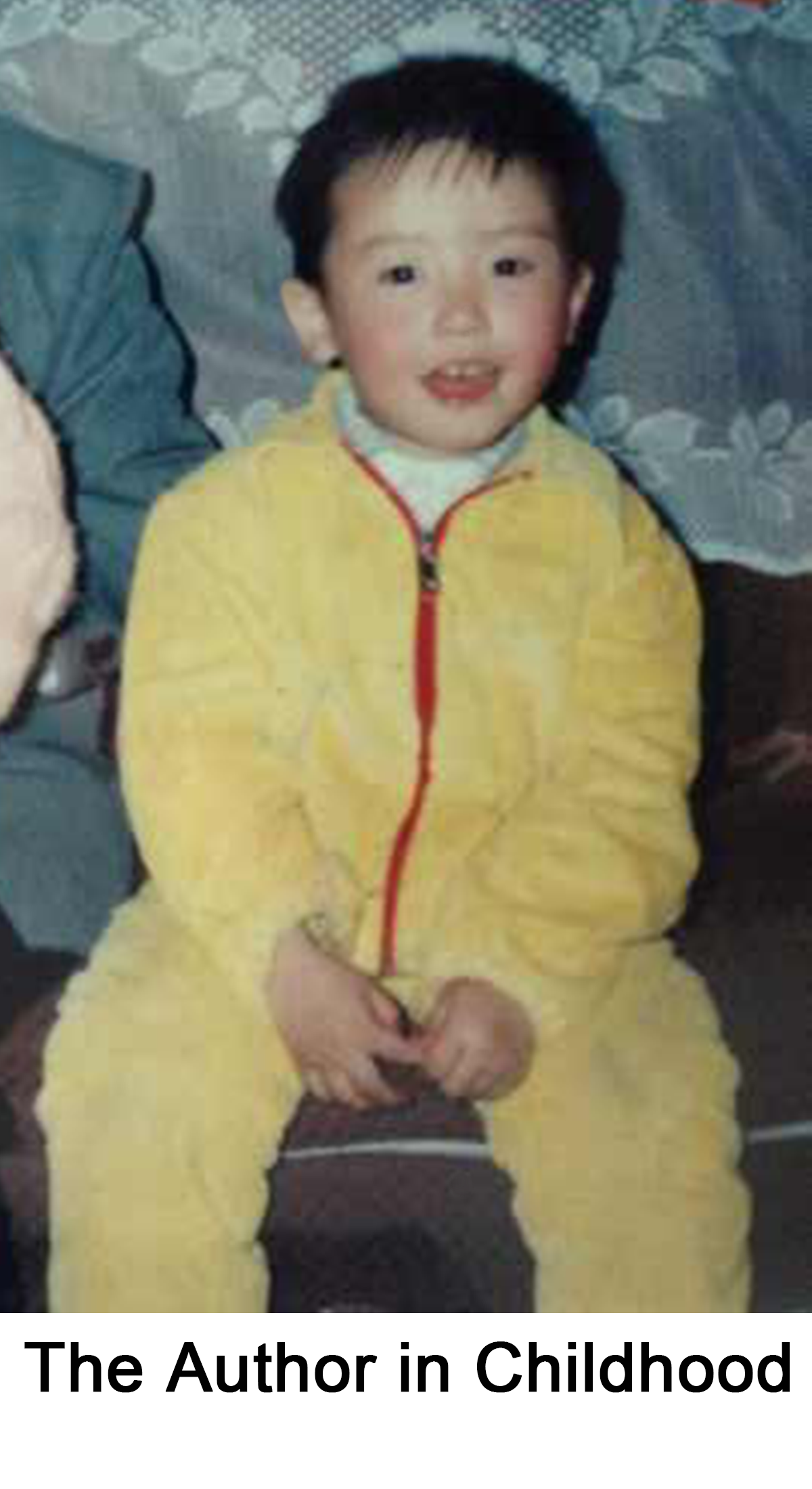}
   \includegraphics[width=12cm]{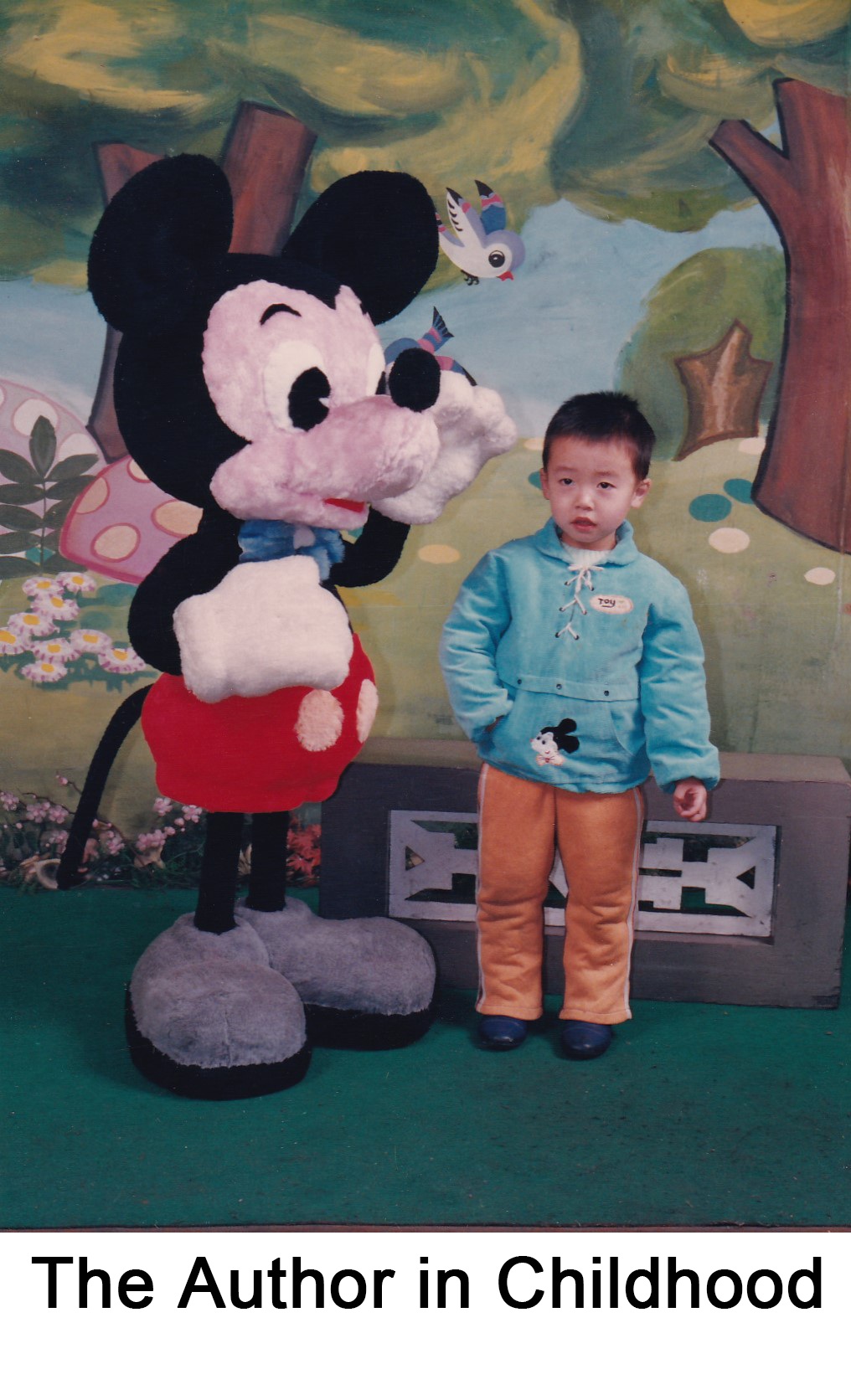}
   \includegraphics[width=12cm]{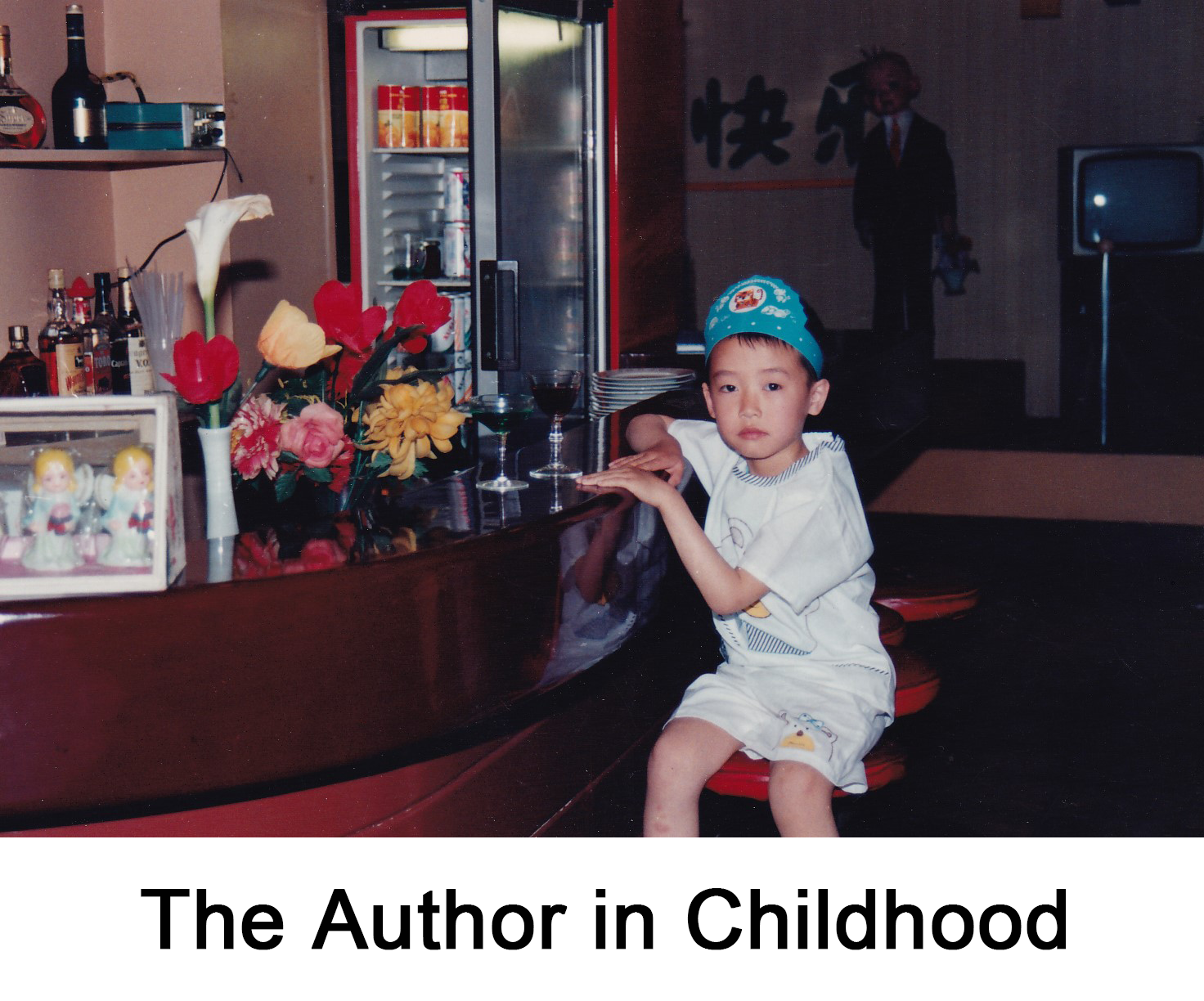}
   \includegraphics[width=12cm]{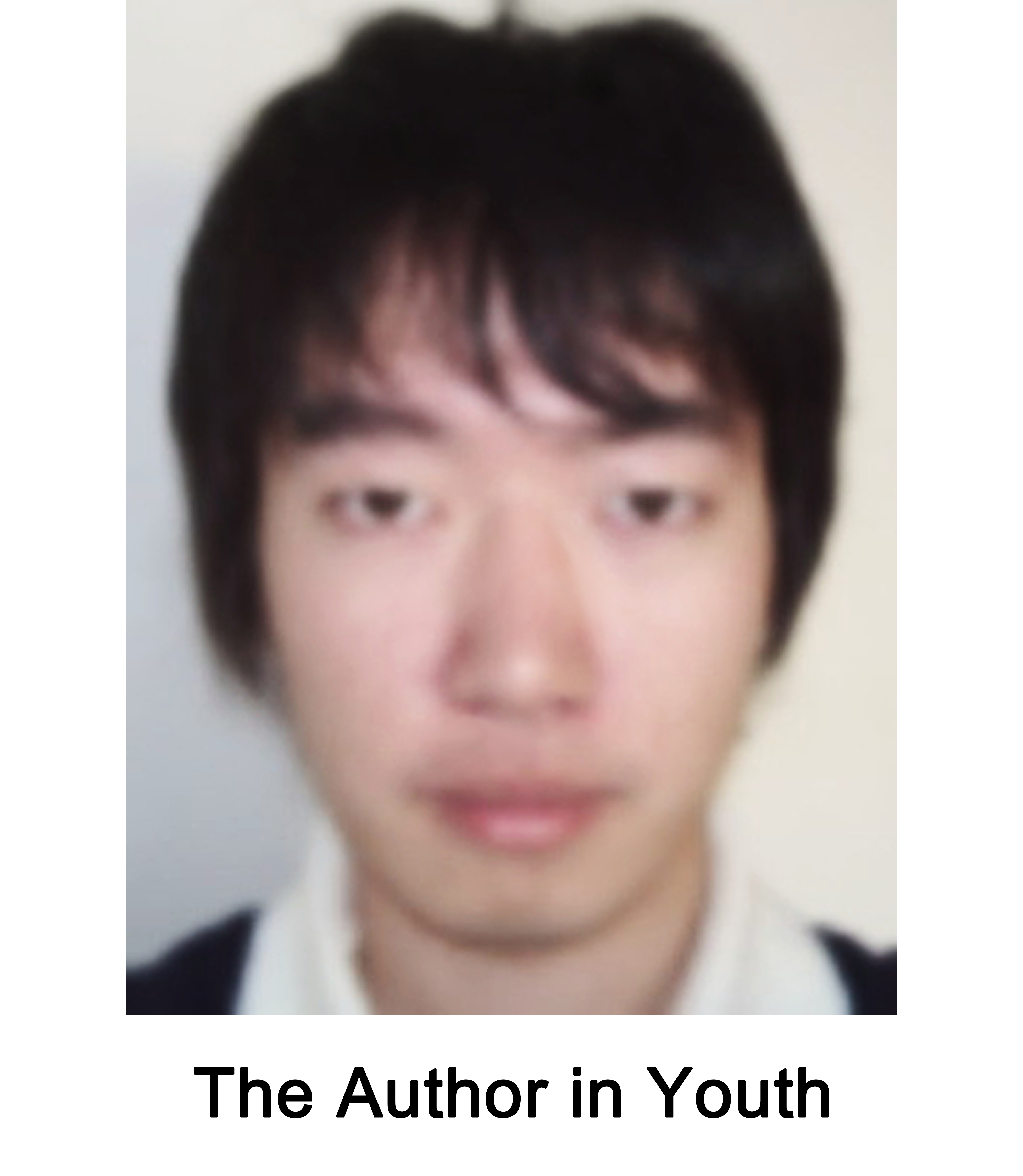}
   \includegraphics[width=12cm]{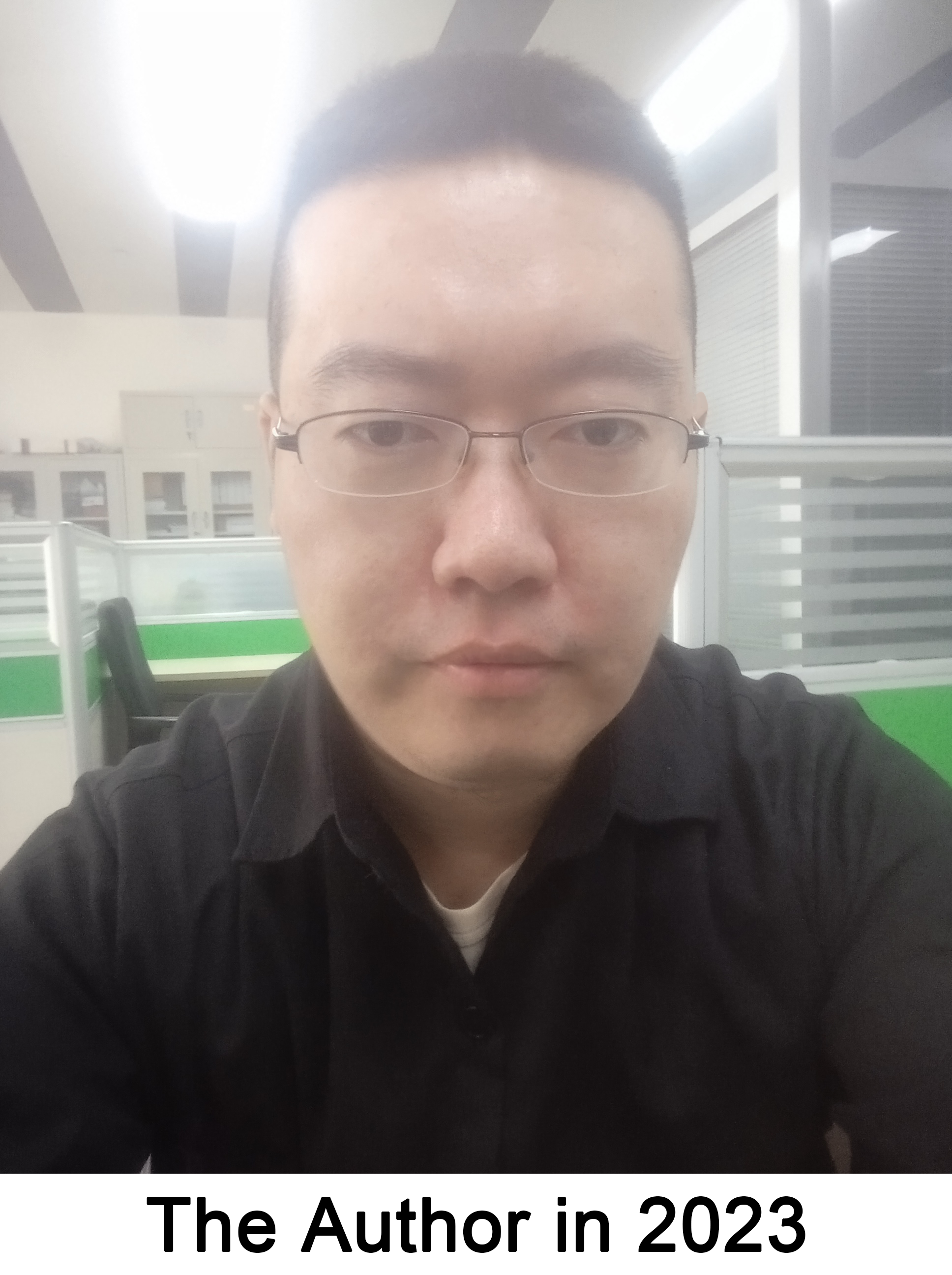}
   \includegraphics[width=12cm]{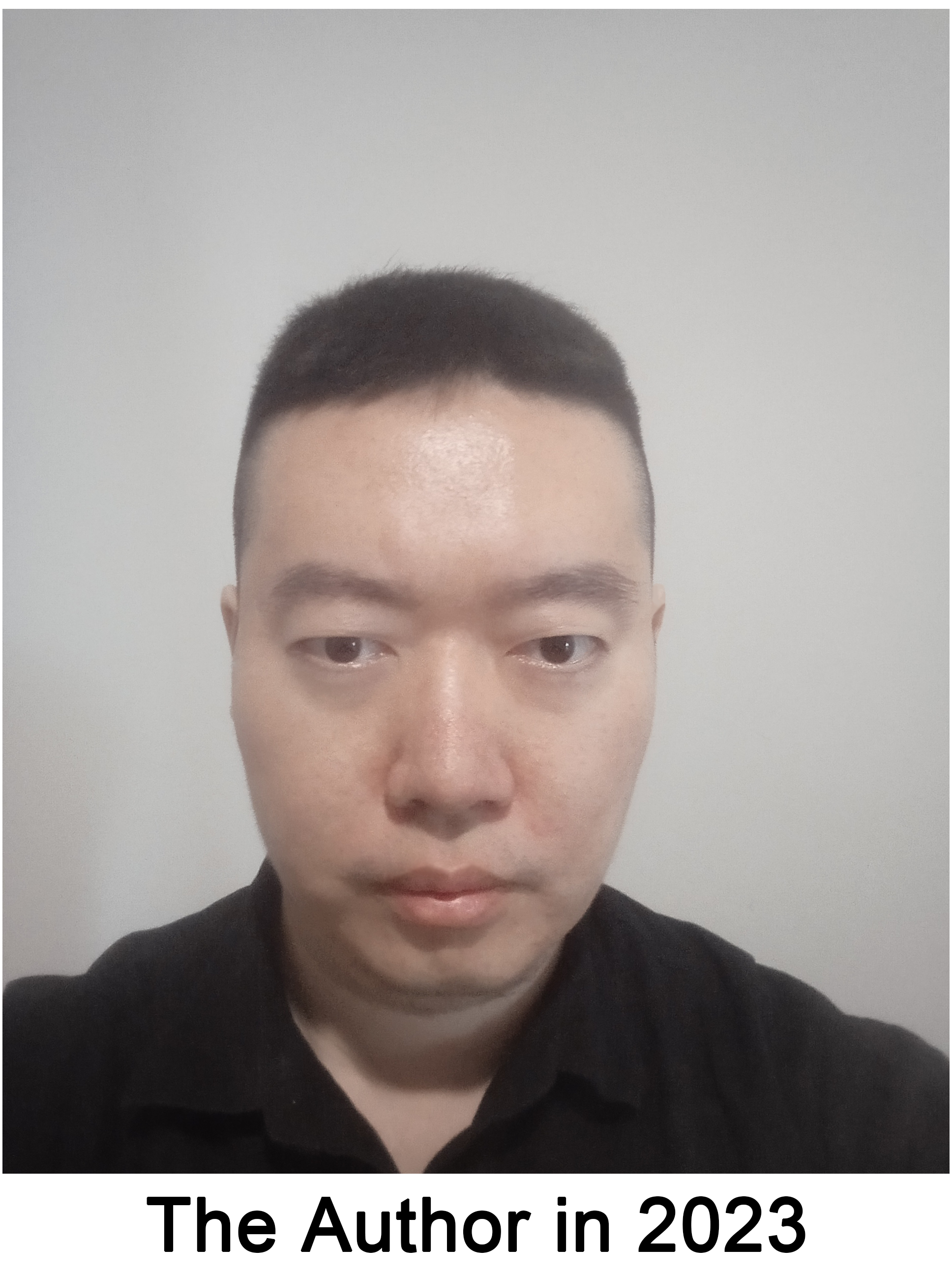}
   \includegraphics[width=12cm]{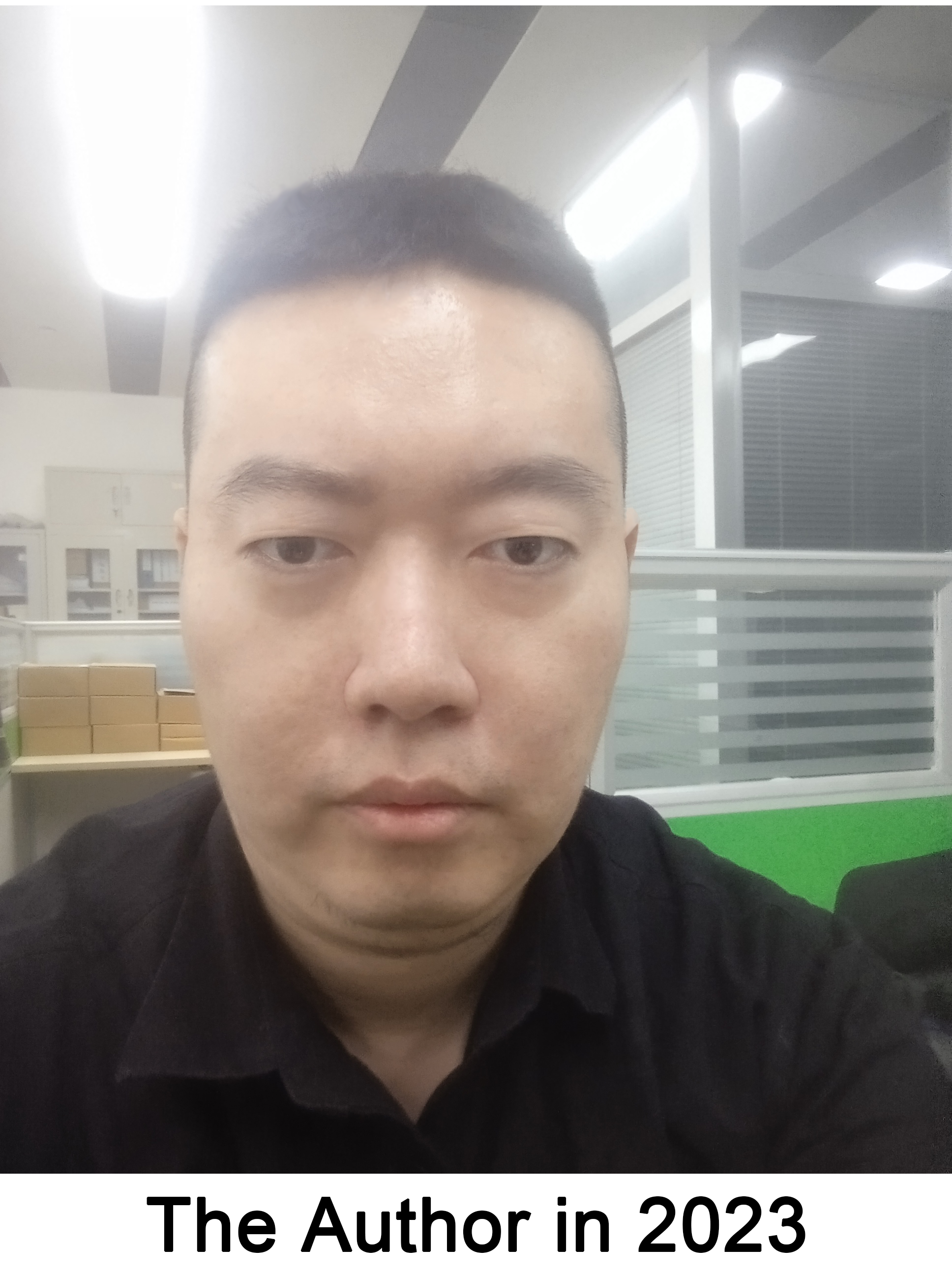}
\end{center}
\end{document}